\documentclass[conference]{IEEEtran}
\usepackage{upgreek}             

\usepackage{ifpdf}

\ifCLASSINFOpdf
   \usepackage[pdftex]{graphicx}
  \graphicspath{{../pdf/}{../jpeg/}}
  \DeclareGraphicsExtensions{.pdf,.jpeg,.png}
\else
  \usepackage[dvips]{graphicx}
 \graphicspath{{../eps/}}
 \DeclareGraphicsExtensions{.eps}
\fi
\usepackage[cmex10]{amsmath}
\begin{document}

%
\title{Complexities Approach to Two Problems In Number Theory}

\author{\IEEEauthorblockN{Bai Yang\\ and Wang Xiuli}
\IEEEauthorblockA{Mount Mercy University\\
Cedar Rapids, IA 52402 USA\\ and Anhui University\\
Hefei, Anhui 230039 China\\
Telephone: (86) 150-560-535-19\\
Email: wangxiuli@ahu.edu.cn}}
\maketitle

\begin{abstract}
By Kolmogorov Complexity,two number-theoretic problems are solved in different way than before,one problem is Maxim Kontsevich and Don Bernard Zagier's Problem 3 \emph{Exhibit  at  least one number  which  does  not  belong  to} $ \mathcal{P}$ (period number) in their paper,another is the problem about existence of bounded coefficients of continued fraction expansion of transcendental number.Thus we show a new approach to mathematical problems in the non-logical discipline.Futhermore,we show that resource-bounded Kolmogorov Complexity and  computational complexity can at least provide tips or principles to  mathematical problems in the non-traditional or logical discipline.
\end{abstract}

%
\IEEEpeerreviewmaketitle

\bibliographystyle{plain}

\section{Introduction}

Maxim Kontsevich and Don Bernard Zagier present a problem ~\cite{Kontsevich2001}, which states:

Problem 3 Exhibit at least one number which does not belong to P(period number).

Later on,They comment:"Each of those problem look very hard and is likely to remain open a long time"


Now we give an example to solve and close the problem by Kolmogorov Complexity and present a more subtle question by Kolmogorov Complexity.

It is easy to show there exist transcendental numbers with regular continued fraction expansion of coefficients   unbounded,but transcendental numbers  with  regular expansion of coefficients bounded are difficult to find.Mailet give first examples ~\cite{Maillet1906},And late on,there are lot of works about such a kind of examples and methods.

In this paper, we give proof for that there exists transcendental numbers of regular continued fraction expansion with bounded coefficients  by Kolmogorov Complexity.

Kolmogorov Complexity is originated from Kolmogorov's research on randomness,Solomonoff's research on induction and Chaitin's research on program length~\cite{MingLiVitanyi1997}.There are variants of definition~\cite{MingLiVitanyi1997},here are several definition of notions relating to this paper and Kolmogorov Complexity.

\emph{Definition} $S$ is set of code word,if $\forall w_i  \in S,\forall w_j \in S$ $w_i$ is not prefix of $w_j$ ,then set of code words is prefix-free code.

\emph{Definition} A prefix-free Turing Machine is Turing Machine which accepts only prefix-free code as program.

\emph{Definition} prefix-free Kolmogorov Complexity of a string, $x$ Given a prefix-free $y$ is $K(x|y)$ which is the length of the shortest prefix-free program that with $y$ as input,outputs $x$.

\emph{Theorem} $K(x|y)$ is non-computable~\cite{MingLiVitanyi1997}.

\emph{Definition} An infinite sequence $r$ is random iff given $x$ any initial segment of the sequence $r$,there is a constant c such that $$K(x)\geq \log |x| +c ,|x| \text{ is the length of $x$}$$

\emph{Theorem} the Lebsgues measure of An infinite sequence over an alphabet $\Sigma =\{a,b\}$ is $1$

\emph{Theorem}  random infinite sequences are non-computable~\cite{MingLiVitanyi1997}.

\section{Computable Number and Non-computable Number}

There are these notions: Turing machine,Computability,computable,partially computable and computable set,computably enumerable set that are relevant to the argument and the example,see ~\cite{Rogers:RecursiveFunctions,Soare:1987:RES:22895} for reference.

\subsection{Computable Reals and computability}
We give the definition of computable reals.

\emph{Definition} if Given $$r \in \mathcal{R}, \forall i\in \mathcal{N}$$,there is Turing machine $\mathbf{M}$ that outputs $i$ bit of $r$,$r$ is computable reals.

This definition can be  extended to complex numbers.

Obviously,Any algebraic numbers including natural numbers and rational numbers are computable.And some transcendental numbers like $\pi,e$ are computable numbers.

Since there are only countable Turing machines,the computable numbers is a countable set.And the set of real numbers is uncountable.Hence we have:

\emph{Theorem} there are uncountable real numbers that are not computable.

More precisely,

\emph{Theorem} there are uncountable transcendental numbers that are not computable.

Computable numbers is related to total computable function or halting Turing machine.Therefore,there has to be an notion that is related to partially computable function.

\emph{Definition} A real number $x = 0.x_1x_2 \dots$ is lower semicomputable if the set of rationals below $x$ is recursively enumerable.
 A number $-x$ is upper semicomputable if $x$ is lower semicomputable. A number $x$ is computable,equivalently, recursive, if it is both lower semicomputable and upper semicomputable ~\cite{MingLiVitanyi1997}

\emph{Theorem}  there are lower semicomputable reals and upper semicomputable reals

\emph{Proof}:All Turing machines can be computably enumerable,let ${T_u(i)},i\in \mathcal{N}$ is a Turing machine that enumerates all  Turing machines.Consider a real $H$ constructed in the following way:
If the $i$ Turing machine halts,then the bit of $H$ in position $i$ is set to be $1$;the bit in position $i$ is set to be $0$,otherwise. Then $H$ is lower semicomputable

\subsection{Chaitin's Number Is Transcendental Number And Not Computable}
For the definition of Chaitin constant or halting probability ,see ~\cite{MingLiVitanyi1997}

\emph{Definition} Chaitin's Number or halting probability is a real number ~\cite{MingLiVitanyi1997}
$$\Omega =\sum_{U(p) < \infty }2^{-p} $$ where $U$ is a prefix-free universal Turing Machine,and $p$ is length of the program that the $U$ halts eventually when inputting it

Obviously It is lower semicomputable and is not computable.

Since all algebraic numbers are computable,Chaitin constant or halting probability and $H$ are not computable,or semicomputable.So,those two numbers are transcendental and non-computable or semicomputable. And all random real number are non-computable.

\section{Period Numbers Are All Computable}

one formula in definition of period numbers is  $$p=\int_{\Delta}\frac{f(x_1,\dots,x_n)}{g(x_1,\dots,x_n)}dx_1\dots dx_n$$
Here $ f$ and $g $ are polynomials with coefficients in
$\mathcal{Q}$, and the integration domain  $\Delta \in \mathcal{R}^n $is given by
polynomial inequalities with rational coefficients(cited from "what is a period")~\cite{Kontsevich2001}.

Obviously,by numerical computation of  integration or Numerical integration,and the definition of formula of  periods,period numbers  are all computable. 

\subsection{Chaitin's Constant Is Not Period Number}
Since Chaitin's Constant is not computable,and is not algebraic number,and all periods are computable,we have

\emph{Theorem} Chaitin's constant is not period number and it is transcendental constant.

Thus close the problem 3 in Kontsevich's paper~\cite{Kontsevich2001}. Moreover, there are non-countable real numbers that are not  periods, since
set of random real number has measure 1 .

\section{There Are Infinite Transcendental Number Of  Simple Continued Fraction With Bounded Coefficients}

\subsection{ Random Sequences of integers map to Bounded Coefficients Of Simple Continued Fraction}
Given a set of finite integers $\mathcal{Z},card(\mathcal{Z}) \geq 2$,the cardinal of  the set of infinite sequences over $\mathcal{Z}$ is  $2^{\aleph}$ and the Lebsgues measure of the set of  random infinite sequences is $1$

Given a random sequence,we can construct a sequence of integers which is not computable or random.Similarly,Given a  bounded random sequence ,we can construct a sequence of  bounded integers which is not computable or random.

\emph{Example} By supremum of Specker sequence, we can construct a real in the following way: if the $n$th digital of the supremum is $0$,make $1$ as the $n$th coefficient of the continued fraction;otherwise,make $2$ as the $n$th coefficient of the continued fraction. The corresponding continued fraction of the supremum is  non-computable 
 In the same way,non-computable or random number will make the corresponding coefficients of the continued fraction non-computable or random. By a random number, we get a non-computable continued fraction with bounded coefficients.
\subsection{There Are Uncountable Transcendental Number   Of Continued Fraction Expansion With Bounded Coefficients}

In the same way,given a random sequence,we can construct a sequence of coefficients of continued fraction which is not computable or random.Similarly,Given a  bounded random sequence ,we can construct a sequence of  bounded coefficients of continued fraction which is not computable or random.

Since all algebraic numbers are computable,those numbers with bounded coefficients of continued fraction which is not computable or random are not computable,thus are all transcendental numbers.We get,

 \emph{Theorem}  There are infinite or uncountable transcendental numbers with bounded coefficients of continued fraction which is not computable or random.
\section{Conclusion and Discussion}
We know such a hierarchy of reals as in the figure 1
\subsection{Refining Some Problems Relating To Transcendental And Algebraic By Complexity}
The proofs about the two problems above  show that it is  very easy to solve those kind of problems by Kolmogorov Complexity,sometimes even trivial.
We can get such results intuitively from the following figure 2:

But when we put computable on the two transcendental problems,they cause  problems like those solved by work of  Mailet~\cite{Maillet1906} and Yosinaga~\cite{Yoshinaga2008},that is,  we have to solve those two existence problems under the computable level,which is much harder than in the whole domain of real number. So  transcendental or algebraic problems, usually become very hard or rich when restricted into the set of computable numbers, and if we have such a problem, we may usually restrict it to the the set of computable numbers or under the computable level to refine it,Figure.2,Figure.3.

Moreover,when we put some restrictions of computational complexity under computability on the kind of problems, those problems almost certainly become very hard,or rich,since computational complexity is  computability time or space of which are restricted.For instance,by intuition and some easy argument ,we can have  such a hierarchy in figure 4 as following:

There are variant computation models over reals, see ~\cite{Weihrauch2000,BlumCuckerShubSmale1998} for reference.here,the computation model and the definition is from ~\cite{Weihrauch2000}.We adopt some result from~\cite{BlumCuckerShubSmale1998}, since the two models are not equivalent,and the model by Blum,Cucker,Shub,Smale is weaker than the Weihrauch's model, and the Weihrauch's model is equivalent to classic model in the computability sense  in the case of $\mathcal{N}$.
Here we know all algebraic numbers and some transcendental numbers are computable in $P$ time, especially from ~\cite{BlumCuckerShubSmale1998},all algebraic numbers  are computable in $P$ time.

If we put real-time computable or linear time computable on the two transcendental problems,they become related to the famous conjecture by Hartmanis and Stearn~\cite{HartmanisSt65}.This can be shown by Figure 5.

Incorporating computational complexity into algebraic and transcendental problem, we can ge such a much more rich hierarchy which motivated a lot of interesting problems. In such a hierarchy,there are a lot of famous problem lied on the boundary among rational,algebraic and transcendental.So,we suggest that any questions relating to rational,algebraic and transcendental,have to incorporate   randomness,computability,and computational complexity and Kolmogorov Complexity,resource-bounded Kolmogorov Complexity.

In addition, we can solve the problem or know some intuitive  approach by randomness, computability and computational complexity. For instance, we can prove Gauss prime theorem by Kolmogorov Complexity. Also we can know from computational complexity, that an algorithm of algebraic number  is in $P$ which implies that another algebraic number by continued fraction is also in $P$

\subsection{Different computation Models}
The two complexity,that is Kolmogorov complexity and computational complexity can solve  problems  by the model in~\cite{Weihrauch2000}: we can refine the input of models as the program $p$,the representation of number $x$ we are computing, length of which is the Kolmogorov complexity of number we are computing, so we can regard the K(x)or $l(p)$ as the $lookahead_M(y,k)$, and the computational complexity is of the algorithm which is the most efficient one. In such a way, the if $\forall y,Lookahead_M(y,k)\rightarrow \infty$ when $n \rightarrow \infty$, the number is non-computable,thus to be a transcendental one, otherwise, it is computable. Under computable, according to $Time_M(y,k)$,can the a hierarchy of reals  be gotten, which intersect or interact with other hierarchy or classification like rational set,algebraic set,transcendental set and reals set, we can refine problems by the complexities and find way or tips to solve them.

Of course, variant models can be applied to number problem, and their results that already exist can pave way or give tips to solve problems

\bibliography{reference}
\end{document}